\documentclass[12pt]{article}
\usepackage{graphicx}

\usepackage{amsmath,amssymb,amstext,amsthm,amsfonts,amsbsy}

\def\be{\begin{eqnarray}}
\def\en{\end{eqnarray}}
\def\nn{\nonumber}
\def\RR{\mathbb R}

\begin{document}

\title{\bf \LARGE A nonlinearly ill-posed problem of reconstructing the temperature from interior data \thanks {Supported by the Council for Natural Sciences of Vietnam} } %
\author{\bf Pham Hoang Quan\thanks{HoChiMinh City National University, Department of Mathematics and Informatics, 227 Nguyen Van Cu, Q. 5,
 HoChiMinh City, VietNam. Email: tquan@pmail.vnn.vn, ddtrong@mathdep.hcmuns.edu.vn}
\addtocounter{footnote}{-1}
- Dang Duc Trong
\footnotemark\\
\bf Alain Pham Ngoc Dinh\thanks{Mathematics Department, Mapmo UMR
6628, BP 67-59, 45067, Orleans Cedex2, France. Email:
alain.pham@math.cnrs.fr}}
\date{}
\maketitle%

{\bf Abstract.} We consider the problem of reconstructing, from the interior data $u(x,1)$, a function $u$ satisfying a nonlinear elliptic equation
\be
\Delta u = f(x,y,u(x,y)),~~~~~x \in \RR, y > 0.
\nn
\en

The problem is ill-posed. Using the method of Green function, the method of Fourier transforms and the method of truncated high frequencies, we shall regularize the problem. Error estimate is given.


{\bf AMS Classification 2000:} 47J06, 35J60, 42A38, 47H10.

\section {Introduction}

In this paper, we consider the problem of reconstructing the temperature of a body from interior measurements. In fact, in many engineering contexts (see, e.g., [BBC]), we cannot attach a temperature sensor at the surface of a body (e.g., the skin of a missile). Hence, to get the temperature distribution on the surface, we have to use the temperature measured inside the body.

Precisely, we consider a two-dimensional body represented by the half-plane $\RR \times \RR^+$. Letting $u(x,y)$ be the temperature of the body at $(x,y) \in \RR\times\RR^+$ and letting $f \equiv f(x,y,u)$ be a (nonlinear) heat source, we have the following nonlinearly nonhomogeneous equation
\be
\Delta u = f, ~~~~~ x \in \RR, y > 0
\label{eq01}
\en
where $\Delta = \frac{\partial^2}{\partial x^2} + \frac{\partial^2}{\partial y^2}$. We assume that the temperature on the line $y = 1$ is known, i. e.,
\be
u(x,1) = \varphi (x),
\label{eq02}
\en
and that
\be
u(x,y) \to 0 ~~\text{as}~~ |x|, y \to \infty.
\label{eq03}
\en

The problem can be referred as a {\it sideways elliptic problem} and the interior measurement $\varphi (x)$ is also called (in geology) the borehole measurement. The problem can be splitted into two problems:

{\bf Problem 1:} finding the function $u$ satisfying
\be
\Delta u = f, ~~~~~ x \in \RR, y > 1
\label{eq04}
\en
subject to the conditions (\ref{eq02}), (\ref{eq03}). Generally, the problem is well-posed. Using the solution of the problem we can calculate the quantity $u_y (x,1) = \phi (x)$.

{\bf Problem 2:} finding a function $u(x,y)~~~~ x \in \RR, 0 < y < 1$\\
satisfying
\be
\Delta u = f, ~~~~~ x \in \RR, 0< y < 1
\label{eq05}
\en
subject to the conditions $u(x,1) =\varphi (x), u_y (x, 1) = \phi (x)$.

The latter problem is a Cauchy elliptic problem and, as known, it is severely ill-posed. Hence, a regularization is in order.

The homogeneous problem ($f\equiv 0$) was studied, by various methods in many papers. Using the mollification method, the homogeneous sideways parabolic problems were considered in [HRS, HR1, AS, L, LV, M] and the references therein. Similarly, the number of papers devoted to the Cauchy problem for linear homogeneous elliptic equation are very rich (see, e.g., [HR, T, B, CHWY, KS]).

Although there are many papers on homogeneous cases, we only find a few papers on nonhomogeneous sideways problems (for both parabolic and elliptic equations). Especially, the papers on the nonlinear case are very rare. In [QT], we use the method of integral equations to consider a sideways elliptic equation with a nonlinear heat source. However, we cannot get an effective method to regularize the problem. In the present paper, we shall consider the problem with a nonlinear heat source $f \equiv f(x,y,u(x,y))$. The remainder of our paper is divided in three sections. In Section 2, using Green functions, we shall transform Problems 1 and 2 into integral equations. In Section 3, we shall prove that Problem 1 has a unique solution. Moreover, we shall give an effective way to approximate the quantity $u_y (x, 1)$. In Section 4, using $\varphi (x), u_y (x,1)$, we shall regularize Problem 2. The main result of the section and of the paper is Theorem 3 (in Subsection 4.2). The method of truncated high frequencies (of Fourier images) will be used and the regularized solution can be found as a fixed point of a contraction. In our knowledge, the latter method is new. Error estimates are given.

\section {Integral equations of Problems 1 and 2}

{\bf 2.1 Problem 1}

Put
\be
N(x,y,\xi,\eta) = -\frac{1}{4\pi} \ln \frac{(x-\xi)^2 +(y-\eta)^2}{(x-\xi)^2 +(y+\eta-2)^2}.
\label{eq06}
\en

For $y > 1, x \in \RR$, integrating the identity
\be
\frac{\partial}{\partial \xi} (u N_\xi - N u_\xi) + \frac{\partial}{\partial \eta} (u N_\eta - N u_\eta) = -N f
\label{eq07}
\en
over the domain $(-m,m) \times (1,n) \backslash B((x,y),\varepsilon)$, where $B((x,y),\varepsilon)$ is the ball with center at $(x,y)$ and radius $\varepsilon > 0$, and letting $n\to \infty, m \to \infty, \varepsilon \to 0$, we get, after some rearrangements,
\be
u(x,y) = A u(x,y),
\label{eq08}
\en
where
\be
A u(x,y) = h(x,y,\varphi) - \int\limits_{-\infty}^{+\infty}\int\limits_{1}^{+\infty}N(x,y;\xi,\eta) f(\xi,\eta,u(\xi,\eta)) d\xi d\eta
\label{eq09}
\en
and
\be
h(x,y,\varphi) = \int\limits_{-\infty}^{+\infty} N_\eta (x,y,\xi,1)\varphi(\xi) d\xi.
\nn
\en

{\bf 2.2 Problem 2}

We repeat $\psi (x) = u_y (x,1)$ with $u$ being an exact solution of Problem 1. We note that we only consider $\varphi$ as an exact data. The function $\psi$ is a processed data dependent on $\varphi$. In fact, in Lemma 1 (Section 3), we shall construct a function $\psi _\varepsilon$ which is an approximation of $\psi$.

Consider
$$\Delta u = f(x,y,u(x,y))\,,\,\,\,\,x \in \RR,\,\,y \in (0,1)$$
subject to the boundary conditions below
\be
u(x,1) = \varphi (x),\,\,\,\,\,x \in \RR,
\label{eq010}
\en
\be
\frac{{\partial u}}{{\partial y}}(x,1) = \psi (x),\,\,\,\,\,x \in \RR.
\nn
\en

We assume in addition that the exact solution $u$ satisfying $\widehat\varphi (\zeta )e^{\left| \zeta  \right|}, \widehat\psi (\zeta )e^{\left| \zeta  \right|}  \in L^2 (\RR)$. We divide this problem into two problems.

{\bf Problem 2.1} Consider the problem
\be
\Delta v = 0, ~~ x \in \RR, ~~ y \in (0,1),
\label{eq011}
\en
\be
v(x,1) = \varphi(x), x \in \RR,
\label{eq012}
\en
\be
\frac{\partial v}{\partial y}(x,1) = \psi (x), x \in \RR.
\label{eq013}
\en

We shall prove that, under appropriate conditions, the problem has a unique solution $v_0$ approximated by a regularized solution $v_\varepsilon$.

{\bf Problem 2.2} Let $u_0$ be an exact solution of $(\ref{eq010})$. If we put $w_0 = u_0 - v_0$ then $w_0$ is the solution of the problem
\be
\Delta w = g(x,y,w),\,\,\,x \in \RR,\,\,\,y \in (0,1)
\label{eq014}
\en
\be
w(x,1) = 0,\,\,\,x \in \RR,
\label{eq015}
\en
\be
\frac{{\partial w}}{{\partial y}}(x,1) = 0,\,\,\,x \in \RR
\label{eq016}
\en
where $g(x,y,w) = f(x,y,w+v_0)$. We shall find a $w_\varepsilon$ which is an approximation of $w_0$ and estimate $\left\| {w_\varepsilon - w_0} \right\|_2$ where $\left\| . \right\|_2$ is the norm in $L^2 (\RR \times (0,1))$. Let $u_\varepsilon = v_\varepsilon + w_\varepsilon$, we shall estimate $\left\| {u_\varepsilon   - u_0 } \right\|_2$.

{\bf 2.2.1 An integral equation of Problem 2.1}

Let
\be
\Gamma (x,y,\xi ,\eta ) =  - \frac{1}{{4\pi }}\ln \left[ {(x - \xi )^2  + (y - \eta )^2 } \right]
\nn
\en
and
\be
G(x,y,\xi ,\eta ) = \Gamma (x,y,\xi ,\eta ) - \Gamma (x, - y,\xi ,\eta ).
\label{eq017}
\en

For $x \in \RR, 0 < y < 1$, integrating the identity
\be
\frac{\partial }{{\partial \xi }}\left( { - vG_\xi   + Gv_\xi  } \right) + \frac{\partial }{{\partial \eta }}\left( { - vG_\eta   + Gv_\eta  } \right) = 0
\nn
\en
over the domain $( - n,n) \times (0,1)\backslash B((x,y),\varepsilon )$ and letting $n \to \infty, \varepsilon \to 0$, we get, after some rearrangements
\be
v(x,y) &=& - \int\limits_{ - \infty }^{ + \infty } {\left[ {\varphi (\xi )G_\eta  (x,y,\xi ,1) - G(x,y,\xi ,1)\psi (\xi )} \right]d\xi } +
\nn\\
&&+ \int\limits_{ - \infty }^{ + \infty } {G_\eta  (x,y,\xi ,0)v(\xi ,0)d\xi}.
\label{eq018}
\en

Letting $y\to 1$ in $(\ref{eq018})$, we have
\be
\frac{1}{\pi }\int\limits_{ - \infty }^{ + \infty } {\frac{1}{{(x - \xi )^2  + 1}}v(\xi ,0)d\xi } + \int\limits_{ - \infty }^{ + \infty } {\left[ { - \varphi (\xi )G_\eta  (x,1,\xi ,1) + G(x,1,\xi ,1)\psi (\xi )} \right]d\xi }  = \varphi (x).
\nn
\en

It can be rewritten as
\be
F_{(1)}  * v_{(0)} (x) = \pi K_{(1)} (x) + \frac{{\sqrt \pi  }}{{\sqrt 2 }}\varphi (x)
\label{eq019}
\en
where
\be
K_{(y)} (x) &=&  - \frac{1}{{\sqrt {2\pi } }}\int\limits_{ - \infty }^{ + \infty } {\left[ { - \varphi (\xi )G_\eta  (x,y,\xi ,1) + G(x,y,\xi ,1)\psi (\xi )} \right]d\xi },
\nn\\
F_{(y)} (x) &\equiv& \frac{y}{{x^2  + y^2 }}, v_{(y)} (x) = v(x,y).
\label{eq020}
\en

Letting
\be
M_{(y,1)} (x) = \frac{{1 - y}}{{x^2  + (y - 1)^2 }} - \frac{{1 + y}}{{x^2  + (y + 1)^2 }}
\nn
\en
and
\be
L_{(\eta ,y)} (x) \equiv \ln \frac{{x^2  + (y - \eta )^2 }}{{x^2  + (y + \eta )^2 }}\,\,\,\,\,(0 < y,\eta  < 1,\,\,x \in \RR),
\nn
\en
we have the Fourier transform of $M, L, F$ as followed
\be
\hat F_{(y)} (\zeta ) &=& \frac{1}{{\sqrt {2\pi } }}\int\limits_{ - \infty }^{ + \infty } {F_{(y)} (x)e^{ - ix\zeta } dx}  = \frac{{\sqrt \pi  }}{{\sqrt 2 }}e^{ - y\left| \zeta  \right|},
\nn\\
\hat L_{(\eta ,y)} (\zeta ) &=& \sqrt {2\pi } \frac{1}{{\left| \zeta  \right|}}\left[ {e^{ - (y + \eta )\left| \zeta  \right|}  - e^{ - \left| {y - \eta } \right|\left| \zeta  \right|} } \right],
\nn\\
\hat M_{(y,1)} (\zeta ) &=& \frac{{\sqrt \pi  }}{{\sqrt 2 }}\left[ {e^{(y - 1)\left| \zeta  \right|}  - e^{ - (y + 1)\left| \zeta  \right|} } \right].
\label{eq021}
\en

From $(\ref{eq020}), (\ref{eq021})$, we get
\be
K_{(y)} (x) =  - \frac{1}{{4\pi }}\left[ {2\varphi  * M_{(y,1)} (x) - \psi  * L_{(1,y)} (x)} \right].
\label{eq022}
\en

From $(\ref{eq019}), (\ref{eq020})$, we have
\be
\hat v_{(0)} (\zeta ) = e^{\left| \zeta  \right|} \left( {\sqrt {2\pi } \hat K_{(1)} (\zeta ) + \hat \varphi (\zeta )} \right).
\label{eq023}
\en

Taking the Fourier transform of (\ref{eq018}), we get
\be
\hat v_{(y)} (\zeta ) &=& e^{\left| \zeta  \right|} \hat F_{(y)} (\zeta )\left( {2\hat K_{(1)} (\zeta ) + \frac{{\sqrt 2 }}{{\sqrt \pi  }}\hat \varphi (\zeta )} \right) - \sqrt {2\pi } \hat K_{(y)} (\zeta )
\nn\\
&=& \frac{1}{2}\hat \varphi (\zeta )\left[ {e^{(1 - y)\left| \zeta  \right|}  + e^{(y - 1)\left| \zeta  \right|} } \right] - \frac{1}{{2\left| \zeta  \right|}}\widehat\psi (\zeta )\left[ {e^{(1 - y)\left| \zeta  \right|}  - e^{(y - 1)\left| \zeta  \right|} } \right] \equiv \aleph (\zeta ,y).
\label{eq024}
\en

{\bf 2.2.2 An integral equation of Problem 2.2}

We recall that $w_0$ an exact solution of Problem $(\ref{eq014})-(\ref{eq016})$.

Let $v_0  \in L^2 (\RR \times (0,1))$ be the exact solution of $(\ref{eq011})-(\ref{eq013})$ and $v_\varepsilon \in L^2 (\RR \times (0,1))$ be a regularized solution.

We write $w_{(y)} (x) = w(x,y)$ and $f_{(\eta ,w_0 ,v_0 )} (\xi ) = f(\xi ,\eta ,v_0 (\xi ,\eta ) + w_0 (\xi ,\eta ))$.

For $x \in \RR, 0 < y < 1$, let $G$ be defined in $(\ref{eq017})$, integrating the identity
\be
\frac{\partial }{{\partial \xi }}( - w_0 G_\xi   + {Gw_0}_\xi  ) + \frac{\partial }{{\partial \eta }}( - w_0 G_\eta   + {Gw_0}_\eta)= Gf_{(\eta ,w_0 ,v_0 )} (\xi )
\label{eq025}
\en
over the domain $( - n,n) \times (0,1)\backslash B((x,y),\varepsilon )$ and letting $n \to\infty, \varepsilon \to 0$, we get, after some rearrangements,
\be
w_0 (x,y) = \int\limits_{ - \infty }^{ + \infty } {w_0 (\xi ,0)G_\eta  (x,y,\xi ,0)d\xi}-\int\limits_{ - \infty }^{ + \infty } {\int\limits_0^1 {G(x,y,\xi ,\eta )f_{(\eta ,w_0 ,v_0 )} (\xi )d\xi d\eta}}.
\label{eq026}
\en

From ($\ref{eq026}$), we have
\be
w_0 (x,y) &=& \frac{1}{\pi }\int\limits_{ - \infty }^{ + \infty } {w_0 (\xi ,0)\frac{y}{{(x - \xi )^2  + y^2 }}d\xi }
\nn\\
&& + \frac{1}{{4\pi }}\int\limits_0^1 {\int\limits_{ - \infty }^{ + \infty } {\ln \frac{{(x - \xi )^2  + (y - \eta )^2 }}{{(x - \xi )^2  + (y + \eta )^2 }}} } f_{(\eta ,w_0 ,v_0 )} (\xi )d\xi d\eta .
\label{eq027}
\en

Letting $y \to 1$, we get
\be
\frac{1}{4}\int\limits_0^1 {\int\limits_{ - \infty }^{ + \infty } {\ln \frac{{(x - \xi )^2  + (1 - \eta )^2 }}{{(x - \xi )^2  + (1 + \eta )^2 }}} } f_{(\eta ,w_0 ,v_0 )} (\xi )d\xi d\eta  + \int\limits_{ - \infty }^{ + \infty } {w_0 (\xi ,0)\frac{1}{{(x - \xi )^2  + 1}}d\xi }  = 0.
\label{eq028}
\en

From $(\ref{eq021})$, we get
\be
&&F_{(y)} (x) \equiv \frac{y}{{x^2  + y^2 }}, L_{(\eta ,y)} (x) \equiv \ln \frac{{x^2  + (y - \eta )^2 }}{{x^2  + (y + \eta )^2 }}\,\,\,\,\,(0 < y,\eta  < 1,\,\,x \in \RR),
\nn\\
&&\hat F_{(y)} (\zeta ) = \frac{{\sqrt \pi  }}{{\sqrt 2 }}e^{ - y\left| \zeta  \right|},
\hat L_{(\eta ,y)} (\zeta ) = \sqrt {2\pi } \frac{1}{{\left| \zeta  \right|}}\left[ {e^{ - (y + \eta )\left| \zeta  \right|}  - e^{ - \left| {y - \eta } \right|\left| \zeta  \right|} } \right].
\label{eq029}
\en

We write
\be
{w_0}_{(y)} (x) \equiv w_0 (x,y).
\nn
\en

From $(\ref{eq028}), (\ref{eq029})$ can be rewritten as
\be
w_{0(0)}*F_{(1)}(x) + \frac{1}{4}\int\limits_0^1 {L_{(\eta ,1)} *f_{(\eta ,w_0 ,v_0 )} (x)d\eta }  = 0.
\nn
\en

Taking the Fourier transform, we have
\be
\widehat{w_{0(0)}} (\zeta ).\hat F_{(1)} (\zeta ) + \frac{1}{4}\int\limits_0^1 {\hat L_{(\eta ,1)} (\zeta )\hat f_{(\eta ,w_0 ,v_0 )} (\zeta )d\eta }  = 0.
\nn
\en

From $(\ref{eq029})$, we have
\be
\widehat{w_{0(0)}} (\zeta ) =  - \frac{1}{2}\int\limits_0^1 {\frac{1}{{\left| \zeta  \right|}}\left[ {e^{ - \eta \left| \zeta  \right|}  - e^{\eta \left| \zeta  \right|} } \right]\hat f_{(\eta ,w_0 ,v_0 )} (\zeta )d\eta}.
\label{eq030}
\en

From $(\ref{eq026})$, we have
\be
w_{0(y)} (x) = \frac{{\sqrt 2 }}{{\sqrt \pi  }}{w_0}_{(0)} *F_{(y)} (x) + \frac{1}{{2\sqrt {2\pi } }}\int\limits_0^1 {L_{(\eta ,y)}  * f_{(\eta ,w_0 ,v_0 )} (x)d\eta }.
\label{eq031}
\en

Taking the Fourier transform of $(\ref{eq031})$, we get
\be
\widehat{w_{0(y)}} (\zeta ) = \frac{{\sqrt 2 }}{{\sqrt \pi  }}\widehat{w_0 }_{(0)} (\zeta ).\hat F_{(y)} (\zeta ) + \frac{1}{{2\sqrt {2\pi } }}\int\limits_0^1 {\hat L_{(\eta ,y)} (\zeta ).\hat f_{(\eta ,w_0 ,v_0 )} (\zeta )d\eta}.
\label{eq032}
\en

From $(\ref{eq029})$ and $(\ref{eq030})$, Eq. $(\ref{eq032})$ takes the form
\be
\widehat{w_0 }_{(y)} (\zeta ) = \frac{1}{2}\int\limits_0^1 {\frac{1}{{\left| \zeta  \right|}}\left[ {e^{(\eta  - y)\left| \zeta  \right|}  - e^{ - \left| {y - \eta } \right|\left| \zeta  \right|} } \right]\hat f_{(\eta ,w_0 ,v_0 )} (\zeta )d\eta}
\label{eq033}
\en
for $\zeta$. We have
\be
\widehat{w_0 }_{(y)} (\zeta ) = \frac{1}{{2\sqrt {2\pi } }}\int\limits_0^1 {\int\limits_{ - \infty }^{ + \infty } {\frac{1}{{\left| \zeta  \right|}}\left[ {e^{(\eta  - y)\left| \zeta  \right|}  - e^{ - \left| {y - \eta } \right|\left| \zeta  \right|} } \right]f_{(\eta ,w_0 ,v_0 )} (\xi )e^{ - i\xi \zeta } d\xi d\eta}}.
\label{eq034}
\en

\section {Finding the solution of Problem 1}

From (\ref{eq08})-(\ref{eq09}), we readily get the following result
\\

{\bf Theorem 1.} {\sl (see [QT]) Suppose that for all $(\xi,\eta,\zeta) \in \RR \times \RR^+ \times \RR$
\be
\left|f'_{\zeta}(\xi,\eta,\zeta)\right| \leq p(\xi,\eta),
\label{eq035}
\en
where $p(\xi,\eta) \in L^1(\RR\times (1,+\infty)), p \geq 0$ satisfies
\be
K \equiv \mathop{\sup}\limits_{(x,y)\in \RR\times (1,+\infty)} \left|\int\limits_{-\infty}^{+\infty}\int\limits_{1}^{+\infty}N(x,y;\xi,\eta) p(\xi,\eta) d\xi d\eta\right| < 1.
\label{eq036}
\en

Put
\be
J=\left\{ u \in C(\RR\times(1,+\infty))\left|\mathop{\lim}\limits_{\sqrt{x^2 + y^2}\to +\infty}u(x,y) = 0\right.\right\}.
\label{eq037}
\en

Then $A: J \to J$ is a contraction and hence $u$ is uniquely determined and can be found by successive approximation.
}
\\

Let $\psi (x) = u_y (x,1)$ with $u$ be an exact solution in the half plane $x \in \RR, y > 1$. From a measured data $\varphi_\varepsilon$ of $\varphi (x) = u(x,1)$, we construct $\psi_\varepsilon$ being an aproximation of $\psi$ and estimate $\left\| {\psi_\varepsilon - \psi} \right\|_{L^2(\RR)}$ by the following lemma
\\

{\bf Lemma 1}

{\sl

Let $u$ be a solution of (\ref{eq08}) as in Theorem 1 and let $\varphi \in L^1(\RR)\cap L^\infty(\RR)$ satisfy \\$\hat\varphi (\zeta) e^{|\zeta|} \in L^2(\RR)$.

Suppose that $f$ satisfies the conditions \be
\left|f(\xi,\eta,\zeta_1) - f(\xi,\eta,\zeta_2)\right| \leq
\left|p(\xi,\eta)\right| \left| \zeta_1 - \zeta_2 \right| ~~~for ~~
all ~~ (\xi,\eta) \in \RR\times \RR^+, \zeta_1,\zeta_2 \in \RR \nn
\en where $p \in L^1(\RR\times (1,+\infty))$,

$$K \equiv \mathop {\sup }\limits_{(x,y) \in \RR \times (1, + \infty )} \int\limits_{ - \infty }^{ + \infty } {\int\limits_1^{ + \infty } {N(x,y;\xi ,\eta )\left| {p(\xi ,\eta )} \right|d\xi } d\eta }  < 1$$
and
$$L = \sqrt {\int\limits_{ - \infty }^{ + \infty } {\left[ {\int\limits_1^{ + \infty } {\int\limits_{ - \infty }^{ + \infty } {\frac{{\eta  - 1}}{{(x - \xi )^2  + (1 - \eta )^2 }}\left| {p(\xi ,\eta )} \right|d\xi d\eta } } } \right]^2 dx} }  <  + \infty.$$

For every $0 < \varepsilon < 1$, we call $\varphi_\varepsilon \in L^2(\RR)$ a measured data such that
$$\left\| {\varphi_\varepsilon - \varphi} \right\|_{L^2(\RR)}  < \varepsilon.$$

From $\varphi_\varepsilon$, we can construct a function $\psi_\varepsilon \in L^2(\RR)$ such that
$$\left\| {\psi _\varepsilon   - \psi } \right\|_{L^2 (\RR)}  < C\varepsilon ^{1/2} $$
where $C$ is independent on $\varepsilon$.
}
\\

{\bf Proof}

Let
\be
k(x,y,u) &=&  - \int\limits_{ - \infty }^{ + \infty } {\int\limits_1^{ + \infty } {N(x,y;\xi ,\eta )f(\xi ,\eta ,u(\xi ,\eta ))d\xi } d\eta }
\nn
\en
and
\be
h(x,y,\varphi) = \int\limits_{-\infty}^{+\infty} N_\eta (x,y,\xi,1) \varphi(\xi) d\xi.
\nn
\en

We have $u(x,y) = h(x,y,\varphi ) + k(x,y,u)$.

Put
$$\widetilde{\varphi _\varepsilon  }(x) = \frac{1}{{\sqrt {2\pi } }}\int\limits_{\left| \zeta  \right| < \frac{1}{{\varepsilon ^{1/2} }}} {\widehat{\varphi _\varepsilon  }(\zeta )e^{i\zeta x} d\zeta }.$$

We have
\be
\widehat{\widetilde{\varphi _\varepsilon  }}(\zeta ) =
\left\{ \begin{array}{l}
 \widehat{\varphi _\varepsilon  }(\zeta )~~~~~~ \left| \zeta  \right| < \frac{1}{{\varepsilon ^{1/2} }} \\
 0  ~~~~~~~~~~~\left| \zeta  \right| \ge \frac{1}{{\varepsilon ^{1/2} }}\, \\
\end{array} \right.
\label{eq038}
\en

Let $u_\varepsilon$ be a solution of $(\ref{eq08})$ with $\varphi$ replaced by $\widetilde{\varphi_\varepsilon}$, i.e.
\be
u_\varepsilon   = h(x,y,\widetilde{\varphi _\varepsilon}) + k(x,y,u_\varepsilon).
\label{eq039}
\en

We denote
\be
h(x,y,\varphi ) &=& h(x,y),
\nn\\
h(x,y,\widetilde{\varphi _\varepsilon}) &=& h_\varepsilon(x,y),
\nn\\
k(x,y,u) &=& k(x,y),
\nn\\
k(x,y,u_\varepsilon  ) &=& k_\varepsilon  (x,y).
\label{eq040}
\en

Let
\be
\psi_\varepsilon  (x) = h_{\varepsilon y} (x,1) + k_{\varepsilon y} (x,1).
\label{eq041}
\en

We have
\be
h(x,y) &=& \int\limits_{ - \infty }^{ + \infty } {N_\eta  (x,y;\xi ,1)\varphi (\xi )d\xi }
\nn\\
&=& \frac{1}{\pi }\int\limits_{ - \infty }^{ + \infty } {\frac{{y - 1}}{{(x - \xi )^2  + (y - 1)^2 }}\varphi (\xi )d\xi}.
\label{eq042}
\en

If we put $F_{(y)} (x) \equiv \frac{y}{{x^2  + y^2 }},y > 0,$ then
\be
\hat F_{(y)} (\zeta ) = \frac{1}{{\sqrt {2\pi } }}\int\limits_{ - \infty }^{ + \infty } {F_{(y)} (x)e^{ - ix\zeta } dx}  = \frac{{\sqrt \pi  }}{{\sqrt 2 }}e^{ - y\left| \zeta  \right|}.
\nn
\en

Taking the Fourier transform of $(\ref{eq042})$, we get
$$\widehat h(\zeta ,y) = \frac{{\sqrt 2 }}{{\sqrt \pi  }}\widehat{F_{(y - 1)} *\varphi }(\zeta ) = \widehat\varphi (\zeta )e^{ - (y - 1)\left| \zeta  \right|} $$
and
$$\widehat h_y (\zeta ,y) =  - \left| \zeta  \right|\widehat\varphi (\zeta )e^{ - (y - 1)\left| \zeta  \right|}.$$

Similarly, we have
\be
\widehat{h_\varepsilon  }(\zeta ,y) = \widehat{\widetilde\varphi }(\zeta )e^{ - (y - 1)\left| \zeta  \right|} .
\nn
\en

We shall find an estimation of $\left\| {h_{\varepsilon y} (.,1) - h_y (.,1)} \right\|_{L^2 (\RR)} $.

Using the inequality $u^4  < e^{2u} ~~~~ \forall u > 1$, we have
\be
\left\| {h_{\varepsilon y} (.,1) - h_y (.,1)} \right\|_{L^2 (\RR)}^2  &=& \int\limits_{ - \infty }^{ + \infty } {\left| \zeta  \right|^2 \left| {\widehat\varphi (\zeta ) - \widehat{\widetilde{\varphi _\varepsilon  }}(\zeta )} \right|^2 d\zeta }
\nn\\
&\le& \int\limits_{\left| \zeta  \right| < \frac{1}{{\varepsilon ^{1/2} }}} {\left| \zeta  \right|^2 \left| {\widehat\varphi (\zeta ) - \widehat{\varphi _\varepsilon  }(\zeta )} \right|^2 d\zeta }  + \int\limits_{\left| \zeta  \right| > \frac{1}{{\varepsilon ^{1/2} }}} {\frac{{e^{2\left| \zeta  \right|} }}{{\left| \zeta  \right|^2 }}\left| {\widehat\varphi (\zeta )} \right|^2 d\zeta}.
\nn
\en

Therefore
$$\left\| {h_{\varepsilon y} (.,1) - h_y (.,1)} \right\|_{L^2 (\RR)}^2  < \varepsilon  + \varepsilon \left\| {e^{\left| \zeta  \right|} \widehat\varphi (\zeta )} \right\|_{L^2 (\RR)}^2 = C_1^2 \varepsilon$$
where $ C_1  = \sqrt {1 + \left\| {e^{\left| \zeta  \right|} \widehat\varphi (\zeta )} \right\|_{L^2 (\RR)}^2 }$.

Hence
\be
\left\| {h_{\varepsilon y} (.,1) - h_y (.,1)} \right\|_{L^2 (\RR)}^{}  < C_1 \sqrt \varepsilon.
\label{eq043}
\en

We have
$$k(x,y) =  - \int\limits_1^{ + \infty } {\int\limits_{ - \infty }^{ + \infty } {N(x,y;\xi ,\eta )f(\xi ,\eta ,u(\xi ,\eta ))d\xi } d\eta }. $$

It follows that
$$k_y (x,1) = \frac{1}{\pi }\int\limits_1^{ + \infty } {\int\limits_{ - \infty }^{ + \infty } {\frac{{1 - \eta }}{{(x - \xi )^2  + (1 - \eta )^2 }}f(\xi ,\eta ,u(\xi ,\eta ))d\xi } d\eta }. $$

Similarly, we get
$$k_{\varepsilon y} (x,1) = \frac{1}{\pi }\int\limits_1^{ + \infty } {\int\limits_{ - \infty }^{ + \infty } {\frac{{1 - \eta }}{{(x - \xi )^2  + (1 - \eta )^2 }}f(\xi ,\eta ,u_\varepsilon  (\xi ,\eta ))d\xi } d\eta }.$$

We have
\be
&&\left\| {k_{\varepsilon y} (.,1) - k_y (.,1)} \right\|_{L^2 (\RR)}=
\nn\\
&& = \frac{1}{\pi }\sqrt {\int\limits_{ - \infty }^{ + \infty } {\left( {\int\limits_1^{ + \infty } {\int\limits_{ - \infty }^{ + \infty } {\frac{{1 - \eta }}{{(x - \xi )^2  + (1 - \eta )^2 }}(f(\xi ,\eta ,u(\xi ,\eta )) - f(\xi ,\eta ,u_\varepsilon  (\xi ,\eta ))d\xi } d\eta } } \right)^2 } dx}
\nn\\
&&\le \frac{1}{\pi }\sqrt {\int\limits_{ - \infty }^{ + \infty } {\left[ {\int\limits_1^{ + \infty } {\int\limits_{ - \infty }^{ + \infty } {\frac{{\eta  - 1}}{{(x - \xi )^2  + (1 - \eta )^2 }}\left| {p(\xi ,\eta )} \right|\left| {u - u_\varepsilon  } \right|d\xi } d\eta } } \right]^2 } dx}
\nn\\
&&\le \frac{1}{\pi }\left\| {u - u_\varepsilon  } \right\|_\infty  L
\label{eq044}
\en
where
$$L = \sqrt {\int\limits_{ - \infty }^{ + \infty } {\left[ {\int\limits_1^{ + \infty } {\int\limits_{ - \infty }^{ + \infty } {\frac{{\eta  - 1}}{{(x - \xi )^2  + (1 - \eta )^2 }}\left| {p(\xi ,\eta )} \right|d\xi d\eta } } } \right]^2 dx} }  <  + \infty.$$

Moreover, we have
\be
\left\| {u - u_\varepsilon  } \right\|_\infty   &\le& \mathop {\sup }\limits_{(x,y) \in \RR \times (1, + \infty )} \left| {\int\limits_{ - \infty }^{ + \infty } {N_\eta  (x,y,\xi ,1)\left[ {\varphi (\xi ) - \widetilde{\varphi _\varepsilon  }(\xi )} \right]d\xi } } \right.
\nn\\
&& \left. { - \int\limits_{ - \infty }^{ + \infty } {\int\limits_1^{ + \infty } {N(x,y,\xi ,\eta )\left( {f(\xi ,\eta ,u(\xi ,\eta )) - f(\xi ,\eta ,u_\varepsilon  (\xi ,\eta ))} \right)d\xi d\eta } } } \right|
\nn
\en
\be
&\le& \mathop {\sup }\limits_{(x,y) \in \RR \times (1 + \infty )} \left[ {\int\limits_{ - \infty }^{ + \infty } {\frac{1}{\pi }\frac{{y - 1}}{{(x - \xi )^2  + (y - 1)^2 }}\left| {\varphi (\xi ) - \widetilde{\varphi _\varepsilon  }(\xi )} \right|d\xi } } \right.
\nn\\
&&\left. { + \int\limits_{ - \infty }^{ + \infty } {\int\limits_1^{ + \infty } {N(x,y,\xi ,\eta )\left| {p(\xi ,\eta )} \right|\left| {u(\xi ,\eta ) - u_\varepsilon  (\xi ,\eta )} \right|d\xi d\eta } } } \right]
\nn\\
&\le& \left\| {\varphi  - \widetilde{\varphi _\varepsilon  }} \right\|_{L^\infty  (\RR)}  + K\left\| {u - u_\varepsilon  } \right\|_\infty
\nn
\en
where
\be
K = \mathop {\sup }\limits_{(x,y) \in \RR \times (1 + \infty )} \int\limits_{ - \infty }^{ + \infty } {\int\limits_1^{ + \infty } {N(x,y,\xi ,\eta )\left| {p(\xi ,\eta )} \right|d\xi d\eta } }  \in (0,1).
\label{eq045}
\en

Hence
\be
\left\| {u - u_\varepsilon  } \right\|_\infty   \le \frac{1}{{1 - K}}\left\| {\varphi  - \widetilde{\varphi _\varepsilon  }} \right\|_{L^\infty  (\RR)}.
\label{eq046}
\en

We get
\be
\left\| {\varphi  - \widetilde{\varphi _\varepsilon  }} \right\|_{L^\infty  (\RR)} &=& \mathop {sup}\limits_{x \in \RR} \left| {\widetilde{\varphi _\varepsilon  }(x) - \varphi (x)} \right|
\nn\\
&=& \mathop {sup}\limits_{x \in \RR} \frac{1}{{\sqrt {2\pi } }}\left| {\int\limits_{\left| \zeta  \right| < \frac{1}{{\varepsilon ^{1/2} }}} {\widehat{\varphi _\varepsilon  }(\zeta )e^{i\zeta x} d\zeta }  - \int\limits_{ - \infty }^{ + \infty } {\widehat\varphi (\zeta )e^{i\zeta x} d\zeta } } \right|
\nn\\
&\le& \frac{1}{{\sqrt {2\pi } }}\left[ {\int\limits_{\left| \zeta  \right| < \frac{1}{{\varepsilon ^{1/2} }}} {\left| {\widehat{\varphi _\varepsilon  }(\zeta ) - \widehat\varphi (\zeta )} \right|d\zeta }  + \int\limits_{\left| \zeta  \right| > \frac{1}{{\varepsilon ^{1/2} }}} {\left| {\widehat\varphi (\zeta )} \right|d\zeta } } \right]
\nn\\
&\le& \frac{1}{{\sqrt {2\pi } }}\left[ {\sqrt {\int\limits_{ - \frac{1}{{\varepsilon ^{1/2} }}}^{\frac{1}{{\varepsilon ^{1/2} }}} {d\zeta } \int\limits_{ - \frac{1}{{\varepsilon ^{1/2} }}}^{\frac{1}{{\varepsilon ^{1/2} }}} {\left| {\widehat{\varphi _\varepsilon  }(\zeta ) - \widehat\varphi (\zeta )} \right|^2 d\zeta } }  + } \right.
\nn\\
&&\left. { + \sqrt {\int\limits_{\left| \zeta  \right| > \frac{1}{{\varepsilon ^{1/2} }}} {e^{2\left| \zeta  \right|} \left| {\widehat\varphi (\zeta )} \right|^2 d\zeta } .\int\limits_{\left| \zeta  \right| > \frac{1}{{\varepsilon ^{1/2} }}} {e^{ - 2\left| \zeta  \right|} d\zeta } } } \right]
\nn\\
&\le& \frac{1}{{\sqrt {2\pi } }}\left[ {\frac{{\sqrt 2 }}{{\varepsilon ^{1/4} }}\left\| {\varphi  - \varphi _\varepsilon  } \right\|_{L^2 (\RR)}  + \sqrt {\frac{2}{3}} \left\| {e^{\left| \zeta  \right|} \widehat\varphi (\zeta )} \right\|_{L^2 (\RR)} \varepsilon ^{3/4} } \right]
\nn\\
&<& C_2 \varepsilon ^{3/4}
\label{eq047}
\en
where
\be
C_2  = \frac{1}{{\sqrt {2\pi } }}\left[ {\sqrt 2  + \sqrt {\frac{2}{3}} \left\| {e^{\left| \zeta  \right|} \widehat\varphi (\zeta )} \right\|_{L^2 (\RR)} } \right].
\nn
\en

Inequalities $(\ref{eq044})-(\ref{eq047})$ give
\be
\left\| {k_{\varepsilon y} (.,1) - k_y (.,1)} \right\|_{L^2 (\RR)}  \le \frac{L}{{\pi (1 - K)}}\left\| {\varphi  - \widetilde\varphi _\varepsilon  } \right\|_\infty   < \frac{{LC_2 }}{{\pi (1 - K)}}\varepsilon ^{{3 \mathord{\left/
 {\vphantom {3 4}} \right.
 \kern-\nulldelimiterspace} 4}} .
\label{eq048}
\en

In view of $(\ref{eq038})-(\ref{eq041}), (\ref{eq043})$ and $(\ref{eq048})$, we have
\be
\left\| {\psi _\varepsilon   - \psi } \right\|_{L^2 (\RR)}  &=& \left\| {u_{\varepsilon y} (.,1) - u_y (.,1)} \right\|_{L^2 (\RR)}
\nn\\
&=& \left\| {h_{\varepsilon y} (.,1) + k_{\varepsilon y} (.,1) - h_y (.,1) - k_y (.,1)} \right\|_{L^2 (\RR)}
\nn\\
&\le& \left\| {h_{\varepsilon y} (.,1) - h_y (.,1)} \right\|_{L^2 (\RR)}  + \left\| {k_{\varepsilon y} (.,1) - k_y (.,1)} \right\|_{L^2 (\RR)}
\nn\\
&<& \frac{{LC_2 }}{{\pi (1 - K)}}\varepsilon ^{{3 \mathord{\left/
 {\vphantom {3 4}} \right.
 \kern-\nulldelimiterspace} 4}}  + C_1 \varepsilon ^{{1 \mathord{\left/
 {\vphantom {1 2}} \right.
 \kern-\nulldelimiterspace} 2}}  \le \left( {\frac{{LC_2 }}{{\pi (1 - K)}} + C_1 } \right)\varepsilon ^{{1 \mathord{\left/
 {\vphantom {1 2}} \right.
 \kern-\nulldelimiterspace} 2}}.
\nn
\en
This complete the proof of Lemma 1.
\\

\section {Regularization of problem}

We shall first regularize the Problem 2.1 (Subsection 4.1). Using the approximated solution $v_\varepsilon$ of Problem 2.1, we shall regularize the Problem 2.2 (Subsection 4.2).

{\bf 4.1 Problem 2.1 (Problem (\ref{eq024}))}

We have the following result

{\bf Theorem 1}

{\sl Let $\varphi  \in L^1 (\RR) \cap L^\infty  (\RR)$. Suppose $\widehat\varphi (\zeta )e^{\left| \zeta  \right|}  \in L^2 (\RR), \widehat\psi (\zeta )e^{\left| \zeta  \right|}  \in L^2 (\RR)$. Then Problem $(\ref{eq011})-(\ref{eq013})$ has a unique solution $v_0  \in L^2 (\RR \times (0,1))$.
}

{\bf Proof}

From $(\ref{eq024})$ and the inequality $\frac{{e^{\left| x \right|}  - 1}}{{\left| x \right|}} \le e^{\left| x \right|}$ we have
\be
\left| {\aleph (\zeta ,y)} \right| \le e^{\left| \zeta  \right|} \left| {\widehat\varphi (\zeta )} \right| + e^{\left| \zeta  \right|} \left| {\widehat\psi (\zeta )} \right| \in L^2 (\RR), \text{~for~all~} 0 \le y \le 1
\nn
\en
where $\aleph$ is as in (\ref{eq024}).

Hence Problem $(\ref{eq011})-(\ref{eq013})$ has a (unique) solution $v_0  \in L^2 (\RR \times (0,1))$.

The proof is completed.

{\bf Theorem 2}

{\sl
Let assumptions in Theorem 1 hold.

For every $0 < \varepsilon  < e^{ - 3}$, let $\varphi_\varepsilon \in L^2 (\RR)$ be measured data such that $\left\| {\varphi _\varepsilon - \varphi } \right\|_{L^2 (\RR)}  < \varepsilon$.

If, in addition, the assumptions of Lemma 1 hold and we have $\left| \zeta  \right|\widehat\varphi (\zeta )e^{\left| \zeta  \right|}  \in L^2 (\RR)$, $\left| \zeta  \right|\widehat\psi (\zeta )e^{\left| \zeta  \right|}  \in L^2 (\RR)$.

Then, from $\varphi_\varepsilon$, we can construct a regularized solution $v_\varepsilon$ such that $\left\| {v_0  - v_\varepsilon  } \right\|_2  < D\left( {\ln \frac{1}{\varepsilon }} \right)^{ - 1}$, where $\left\| . \right\|_2$ is the norm in $L^2 (\RR \times (0,1))$ and $D$ is a positive constant independent of $\varepsilon$.}

{\bf Proof}

We recall that $\widetilde{\varphi _\varepsilon}$ is defined in Lemma 1, we get
\be
\left\| {\widetilde{\varphi _\varepsilon  } - \varphi } \right\|_{L^2 (\RR)}^2  &=& \int\limits_{\left| \zeta  \right| < \varepsilon ^{ - 1/2} } {\left| {\widehat{\varphi _\varepsilon  }(\zeta ) - \widehat\varphi (\zeta )} \right|^2 d\zeta }  + \int\limits_{\left| \zeta  \right| > \varepsilon ^{ - 1/2} } {\left| {\widehat\varphi (\zeta )} \right|^2 d\zeta }
\nn\\
&<& \varepsilon ^2  + \int\limits_{\left| \zeta  \right| > \varepsilon ^{ - 1/2} } {\frac{{\left| {\widehat\varphi (\zeta )} \right|^2 e^{2\left| \zeta  \right|} }}{{\left| \zeta  \right|^4 }}d\zeta }
\nn\\
&<& \varepsilon ^2  + \varepsilon ^2 \left\| {\widehat\varphi (\zeta )e^{\left| \zeta  \right|} } \right\|_{L^2 (\RR)}^2 = \varepsilon ^2 C_3^2
\nn
\en
where $C_3  = \sqrt {1 + \left\| {\widehat\varphi (\zeta )e^{\left| \zeta  \right|} } \right\|_{L^2 (\RR)}^2 }$.

Therefore $\left\| {\widetilde{\varphi _\varepsilon  } - \varphi } \right\|_{L^2 (\RR)}  < C_3 \varepsilon$.

Put
\be
\phi _\varepsilon  (x) = \frac{1}{{\sqrt {2\pi } }}\int\limits_{ - \frac{1}{6}\ln \frac{1}{\varepsilon }}^{\frac{1}{6}\ln \frac{1}{\varepsilon }} {\widehat{\widetilde\varphi }(\zeta )e^{i\zeta x} d\zeta }.
\nn
\en

We have
\be
\widehat{\phi _\varepsilon  }(\zeta ) = \left\{ \begin{array}{l}
 \widehat{\widetilde{\varphi _\varepsilon  }}(\zeta ) ~~~~ \left| \zeta  \right| < \frac{1}{6}\ln \frac{1}{\varepsilon } \\
 0 ~~~~~~~~~ \left| \zeta  \right| \ge \frac{1}{6}\ln \frac{1}{\varepsilon } \\
 \end{array} \right..
\nn
\en

It follows that
\be
\left\| {\left( {\widehat{\phi _\varepsilon  }(\zeta ) - \widehat\varphi (\zeta )} \right)e^{\left| \zeta  \right|} } \right\|_{L^2 (\RR)}^2  &=& \int\limits_{ - \frac{1}{6}\ln \frac{1}{\varepsilon }}^{\frac{1}{6}\ln \frac{1}{\varepsilon }} {\left| {\widehat{\widetilde{\varphi _\varepsilon  }}(\zeta ) - \widehat\varphi (\zeta )} \right|^2 e^{2\left| \zeta  \right|} d\zeta }  + \int\limits_{\left| \zeta  \right| > \frac{1}{6}\ln \frac{1}{\varepsilon }} {\left| {\widehat\varphi (\zeta )} \right|^2 e^{2\left| \zeta  \right|} d\zeta }
\nn\\
&\le& e^{\frac{1}{3}\ln \frac{1}{\varepsilon }} \left\| {\widetilde{\varphi _\varepsilon  } - \varphi } \right\|_{L^2 (\RR)}^2  + \frac{1}{{\left( {\frac{1}{6}\ln \frac{1}{\varepsilon }} \right)^2 }}\left\| {\zeta \widehat\varphi (\zeta )e^{\left| \zeta  \right|} } \right\|_{L^2 (\RR)}^2
\nn\\
&<& C_3^2 e^{\frac{1}{3}\ln \frac{1}{\varepsilon }} \varepsilon ^2  + \frac{1}{{\left( {\frac{1}{6}\ln \frac{1}{\varepsilon }} \right)^2 }}\left\| {\zeta \widehat\varphi (\zeta )e^{\left| \zeta  \right|} } \right\|_{L^2 (\RR)}^2
\nn\\
&<& C_3^2 \varepsilon ^{5/3}  + \frac{{36}}{{\ln ^2 \frac{1}{\varepsilon }}}\left\| {\zeta \widehat\varphi (\zeta )e^{\left| \zeta  \right|} } \right\|_{L^2 (\RR)}^2.
\nn
\en

Using the inequality $\varepsilon ^{5/3}  < \frac{1}{{\ln ^2 \frac{1}{\varepsilon }}}$ as $\varepsilon  < e^{-3}$, we get
\be
\left\| {\left( {\widehat{\phi _\varepsilon  }(\zeta ) - \widehat\varphi (\zeta )} \right)e^{\left| \zeta  \right|} } \right\|_{L^2 (\RR)}^2  < C_4^2 \frac{1}{{\ln ^2 \frac{1}{\varepsilon }}},
\nn
\en
where
\be
C_4  = \sqrt {C_3^2  + 36\left\| {\zeta \widehat\varphi (\zeta )e^{\left| \zeta  \right|} } \right\|_{L^2 (\RR)}^2 }.
\label{eq049}
\en

Using lemma 1, there exists $\psi _\varepsilon \in L^2 (\RR)$ such that
\be
\left\| {\psi _\varepsilon   - \psi } \right\|_{L^2 (\RR)}  < C\varepsilon ^{1/2}.
\nn
\en

Put
\be
\Psi _\varepsilon  (x) = \frac{1}{{\sqrt {2\pi } }}\int\limits_{ - \frac{1}{6}\ln \frac{1}{\varepsilon }}^{\frac{1}{6}\ln \frac{1}{\varepsilon }} {\widehat{\Psi _\varepsilon  }(\zeta )e^{i\zeta x} d\zeta}.
\nn
\en

We have
\be
\widehat{\Psi _\varepsilon  }(\zeta ) = \left\{ \begin{array}{l}
 \widehat{\psi _\varepsilon  }(\zeta ) ~~~~ \left| \zeta  \right| < \frac{1}{6}\ln \frac{1}{\varepsilon } \\
 0 ~~~~~~~~~\left| \zeta  \right| \ge \frac{1}{6}\ln \frac{1}{\varepsilon } \\
 \end{array} \right..
\nn
\en

Similarly, using the inequality $\varepsilon ^{2/3}  < \frac{1}{{\ln ^2 \frac{1}{\varepsilon}}}$ as $\varepsilon < e^{-3}$, we get
\be
\left\| {\left( {\widehat{\Psi _\varepsilon  }(\zeta ) - \widehat\psi (\zeta )} \right)e^{\left| \zeta  \right|} } \right\|_{L^2 (\RR)}^2  < C_5^2 \frac{1}{{\ln ^2 \frac{1}{\varepsilon}}},
\nn
\en
where
\be
C_5  = \sqrt {C^2  + 36\left\| {\zeta \widehat\psi (\zeta )e^{\left| \zeta  \right|} } \right\|_{L^2 (\RR)}^2 }.
\label{eq050}
\en

We put \be \aleph _\varepsilon  (\zeta ,y) =
\frac{1}{2}\widehat{\phi _\varepsilon  }(\zeta )\left[ {e^{(1 -
y)\left| \zeta  \right|}  + e^{(y - 1)\left| \zeta  \right|} }
\right] - \frac{1}{{2\left| \zeta  \right|}}\widehat{\Psi
_\varepsilon  }(\zeta )\left[ {e^{(1 - y)\left| \zeta  \right|}  -
e^{(y - 1)\left| \zeta  \right|} } \right] \nn \en and \be
v_\varepsilon  (x,y) = \frac{1}{{\sqrt {2\pi } }}\int\limits_{ -
\infty }^{ + \infty } {\aleph _\varepsilon  (\zeta ,y)e^{i\zeta x}
d\zeta}. \label{eq051} \en

From $(\ref{eq024}), (\ref{eq051})$, we get \be &&\left\|
{v_\varepsilon   - v_0 } \right\|_2  = \left\| {\aleph _\varepsilon
- \aleph } \right\|_2
\nn\\
&&~~\le \frac{1}{2}\left( {\left\| {\left( {\widehat\varphi (\zeta ) - \widehat{\phi _\varepsilon  }(\zeta )} \right)\left[ {e^{(1 - y)\left| \zeta  \right|}  + e^{(y - 1)\left| \zeta  \right|} } \right]} \right\|_2 + \left\| {\frac{{\left( {\widehat\psi (\zeta ) - \widehat{\Psi _\varepsilon  }(\zeta )} \right)}}{{\left| \zeta  \right|}}\left[ {e^{(1 - y)\left| \zeta  \right|}  - e^{(y - 1)\left| \zeta  \right|} } \right]} \right\|_2 } \right)
\nn\\
&&~~\le \left\| {(\widehat\varphi (\zeta ) - \widehat{\phi _\varepsilon  }(\zeta ))e^{\left| \zeta  \right|} } \right\|_2  + \left\| {(\widehat\psi (\zeta ) - \widehat{\Psi _\varepsilon  }(\zeta ))e^{\left| \zeta  \right|} } \right\|_2
\nn\\
&&~~< (C_4  + C_5 )\left( {\ln \frac{1}{\varepsilon }} \right)^{ - 1}.
\nn
\en

Hence $\left\| {v_0  - v_\varepsilon  } \right\|_2  < D\left( {\ln \frac{1}{\varepsilon }} \right)^{-1}$.

The proof is completed.

{\bf 4.2 Problem 2.2 (Problem (\ref{eq034}))}

{\bf Theorem 3}

{\sl
Let assumptions in Theorem 1 and Theorem 2 hold.

Let $v_0 \in L^2 (\RR \times (0,1))$ be an exact solution of $(\ref{eq011})-(\ref{eq013})$ and $v_\varepsilon   \in L^2 (\RR \times (0,1))$ be the regularized solution of $(\ref{eq011})-(\ref{eq013})$.

Suppose that $f$ satisfies the conditions $(\ref{eq027})$ and
$\left| {p(\xi ,\eta )} \right| \le k\,\,\forall (\xi ,\eta ) \in \RR \times (0,1)$.

Assume in addition that the exact solution $u_0 = v_0  + w_0\in L^2 (\RR \times (0,1))$ of $(\ref{eq010})$ satisfying
\be
e^{3\left| \zeta  \right|} \left| {\widehat f_{(\eta ,v_0 ,w_0 )} (\zeta )} \right| \in L^2 (\RR \times (0,1)).
\nn
\en

Then there exists a regularized solutions $w_\varepsilon$ and $u_\varepsilon = v_\varepsilon + w_\varepsilon$ of $(\ref{eq010})$ such that
\be
\left\| {w_\varepsilon   - w_0 } \right\|_2  < C\left( {\ln \frac{1}{\varepsilon }} \right)^{-1/2}
\nn
\en
and
\be
\left\| {u_\varepsilon   - u_0 } \right\|_2  < E\left( {\ln \frac{1}{\varepsilon }} \right)^{-1/2}
\nn
\en
where $C$ and $E$ independent of $\varepsilon$.
}

{\bf Proof}

From $(\ref{eq033})$, we have
\be
\widehat{w_0 }_{(y)} (\zeta ) &=& \frac{1}{2}\int\limits_0^1 {\frac{1}{{\left| \zeta  \right|}}\left[ {e^{(\eta  - y)\left| \zeta  \right|}  - e^{ - \left| {y - \eta } \right|\left| \zeta  \right|} } \right]\widehat f_{(\eta ,w_0 ,v_0 )} (\zeta )d\eta }
\nn\\
& =& \frac{1}{2}\int\limits_y^1 {\frac{1}{{\left| \zeta  \right|}}\left[ {e^{(\eta  - y)\left| \zeta  \right|}  - e^{(y - \eta )\left| \zeta  \right|} } \right]\widehat f_{(\eta ,w_0 ,v_0 )} (\zeta )d\eta }.
\nn
\en

We put \be T\left( {w^{(\alpha )} (x,y)} \right) = \frac{1}{{2\sqrt
{2\pi } }}\int\limits_y^1 {\int\limits_{ - \alpha }^\alpha
{\frac{1}{{\left| \zeta  \right|}}\left[ {e^{(\eta  - y)\left| \zeta
\right|}  - e^{(y - \eta )\left| \zeta  \right|} } \right]\widehat
f_{(\eta ,w^{(\alpha )} ,v_\varepsilon  )} (\zeta )e^{i\zeta x}
d\zeta } d\eta } \label{eq052} \en

We shall prove \be \left\| {T^m (w^{(\alpha )} (.,y)) - T^m
(w^{1(\alpha )} (.,y))} \right\|_{L^2 (\RR)}^2  \le \left( {k^2
e^{2\alpha } } \right)^m \frac{{(1 - y)^m }}{{m!}}\left| {\left\|
{w^{(\alpha )}  - w^{1(\alpha )} } \right\|} \right|^2
\label{eq053}
\en
 for every $\alpha  > 0, y \in (0,1), m \ge 1$ and $w^{(\alpha )}
,w^{1(\alpha )} \in C([0,1];L^2 (\RR))$, where $| \| . \| |$ denotes
sup-norm in $C([0,1];L^2 (\RR))$.

We shall prove $(\ref{eq051})$ by induction.

We define
\be
\chi_{[ - \alpha ,\alpha ]} (\zeta ) = \left\{ \begin{array}{l}
 1 ~~~~~ |\zeta| \leq \alpha\\
 0 ~~~~~ |\zeta| > \alpha\\
 \end{array} \right..
\nn
\en

For $m = 1$, noting that
\be
\widehat{Tw_{}^{(\alpha )} }(\zeta ,y) = \frac{1}{2}\chi _{[ - \alpha ,\alpha ]} (\zeta )\int\limits_y^1 {\frac{1}{{\left| \zeta  \right|}}\left[ {e^{(\eta  - y)\left| \zeta  \right|}  - e^{(y - \eta )\left| \zeta  \right|} } \right]\widehat f_{(\eta ,w^{(\alpha )} ,v_\varepsilon  )} (\zeta )d\eta },
\nn
\en
we have
\be
&&\left\| {Tw^{(\alpha )} (.,y) - T{w^1}^{(\alpha )} (.,y)} \right\|_{L^2 (\RR)}^2 =
\nn\\
&&~~= \left\| {\widehat{Tw^{(\alpha )} }(.,y) - \widehat{T{w^1}^{(\alpha )} }(.,y)} \right\|_{L^2 (\RR)}^2
\nn\\
&&~~= \frac{1}{4}\left\| {\int\limits_y^1 {\chi _{[ - \alpha ,\alpha ]} (\zeta )\frac{{e^{(\eta  - y)\left| \zeta  \right|}  - e^{(y - \eta )\left| \zeta  \right|} }}{{\left| \zeta  \right|}}\left( {\widehat f_{(\eta ,w^{(\alpha )} ,v_\varepsilon  )} (\zeta ) - \widehat f_{(\eta ,w^{1(\alpha )} ,v_\varepsilon  )} (\zeta )} \right)d\eta } } \right\|_{L^2 (\RR)}^2
\nn\\
&&~~\le \int\limits_{ - \infty }^{ + \infty } {\left| {\int\limits_y^1 {\chi _{[ - \alpha ,\alpha ]} (\zeta )e^{(\eta  - y)\left| \zeta  \right|} \left| {\widehat f_{(\eta ,w^{(\alpha )} ,v_\varepsilon  )} (\zeta ) - \widehat f_{(\eta ,w^{1(\alpha )} ,v_\varepsilon  )} (\zeta )} \right|d\eta } } \right|^2 d\zeta }
\nn
\en
\be
&&~~\le \int\limits_{ - \infty }^{ + \infty } {\left| {\int\limits_y^1 {e^{(\eta  - y)\alpha } \left| {\widehat f_{(\eta ,w^{(\alpha )} ,v_\varepsilon  )} (\zeta ) - \widehat f_{(\eta ,w^{1(\alpha )} ,v_\varepsilon  )} (\zeta )} \right|d\eta } } \right|^2 d\zeta }
\nn\\
&&~~\le e^{2\alpha } (1 - y)\int\limits_{ - \infty }^{ + \infty } {\int\limits_y^1 {\left| {\widehat f_{(\eta ,w^{(\alpha )} ,v_\varepsilon  )} (\zeta ) - \widehat f_{(\eta ,w^{1(\alpha )} ,v_\varepsilon  )} (\zeta )} \right|^2 d\eta } d\zeta }
\nn\\
&&~~\le k^2 e^{2\alpha } (1 - y)\int\limits_y^1 {\int\limits_{ - \infty }^{ + \infty } {\left| {w^{(\alpha )} (x,\eta ) - w^{1(\alpha )} (x,\eta )} \right|^2 dx} d\eta }
\nn\\
&&~~\le e^{2\alpha } (1 - y)k^2 \int\limits_y^1 {\left\| {w^{(\alpha )} (.,\eta ) - w^{1(\alpha )} (.,\eta )} \right\|_{L^2 (\RR)}^2 d\eta }
\nn\\
&&~~\le e^{2\alpha } (1 - y)k^2 \left| {\left\| {w^{(\alpha )}  - w^{1(\alpha )} } \right\|} \right|^2.
\nn
\en

Therefore $(\ref{eq053})$ holds.

Suppose that $(\ref{eq053})$ holds for $m = j$. We shall prove that
$(\ref{eq051})$ holds for $m = j + 1$. We have \be &&\left\| {T^{j +
1} \left( {w^{(\alpha )} (.,y)} \right) - T^{j + 1} \left(
{w^{1(\alpha )} (.,y)} \right)} \right\|_{L^2 (\RR)}^2 =
\nn\\
&&~~= \left\| {\widehat T\left( {T^j \left( {w^{(\alpha )} (.,y)} \right)} \right) - \widehat T\left( {T^j \left( {w^{1(\alpha )} (.,y)} \right)} \right)} \right\|_{L^2 (\RR)}^2
\nn\\
&&~~\le k^2 e^{2\alpha } (1 - y)\int\limits_y^1 {\left\| {T^j \left( {w^{(\alpha )} (.,\eta )} \right) - T^j \left( {w^{1(\alpha )} (.,\eta )} \right)} \right\|_{L^2 (\RR)}^2 d\eta }
\nn\\
&&~~\le k^2 e^{2\alpha } (1 - y)\int\limits_y^1 {\left( {k^2 e^{2\alpha } } \right)^j \frac{{(1 - \eta )^j }}{{j!}}\left| {\left\| {w^{(\alpha )}  - w^{1(\alpha )} } \right\|} \right|^2 d\eta }
\nn
\en
\be
&&~~\le \left( {k^2 e^{2\alpha } } \right)^{j + 1} \frac{{(1 - y)^{j + 1} }}{{(j + 1)!}}\left| {\left\| {w^{(\alpha )}  - w^{1(\alpha )} } \right\|} \right|^2.
\nn
\en

Therefore, by the induction principle, $(\ref{eq053})$ holds for
every $m$. From $(\ref{eq053})$, we get \be \left| {\left\| {T^m
\left( {w^{(\alpha )} } \right) - T^m \left( {w^{1(\alpha )} }
\right)} \right\|} \right|^2 \le \left( {k^2 e^{2\alpha } }
\right)^m \frac{1}{{m!}}\left| {\left\| {w^{(\alpha )}  -
w^{1(\alpha )} } \right\|} \right|^2. \nn \en

Since $\mathop {\lim }\limits_{m \to \infty } \left( {k^2 e^{2\alpha } } \right)^m \frac{1}{{m!}} = 0$, there exists a positive integer number $m_0$ such that
\be
\left( {k^2 e^{2\alpha } } \right)^{m_0 } \frac{1}{{m_0 !}} < 1.
\nn
\en

Hence $T^{m_0}$ is a contraction in $C\left( {[0,1];L^2 (\RR)} \right)$. It follows that the equation $T^{m_0 } (w^{(\alpha )} ) = w^{(\alpha )}$ has a unique solution $w^{(\alpha )}  \in C\left( {[0,1];L^2 (\RR)} \right)$.

We claim that $T(w^{(\alpha )} ) = w^{(\alpha )}$. In fact, one has $T(T^{m_0 } (w^{(\alpha )} )) = T(w^{(\alpha )} )$. Hence $T^{m_0 } (T(w^{(\alpha )} )) = T(w^{(\alpha )})$. By the uniqueness of the fixed point of $T^{m_0 }$, one has $T(w^{(\alpha )} ) = w^{(\alpha )}$, i.e., the equation $T(w^{(\alpha )} ) = w^{(\alpha )}$ has a unique solution $w^{(\alpha )}  \in C\left( {[0,1];L^2 (\RR)} \right)$.

We have
\be
&&4\left\| {w_{0(y)}  - w_{}^{(\alpha )} (.,y)} \right\|_{L^2 (\RR)}^2
\nn\\
&&~~= \int\limits_{ - \alpha }^\alpha  {\left| {\int\limits_y^1 {\frac{{[e^{(\eta  - y)\left| \zeta  \right|}  - e^{(y - \eta )\left| \zeta  \right|} ]}}{{\left| \zeta  \right|}}\left[ {\widehat f_{(\eta ,w_0 ,v_0 )} (\zeta ) - \widehat f_{(\eta ,w^{(\alpha )} ,v_\varepsilon  )} (\zeta )} \right]} d\eta } \right|} ^2 d\zeta
\nn\\
&&~~~~+ \int\limits_{\left| \zeta  \right| > \alpha }^{} {\left| {\int\limits_y^1 {\frac{{[e^{(\eta  - y)\left| \zeta  \right|}  - e^{(y - \eta )\left| \zeta  \right|} ]}}{{\left| \zeta  \right|}}\widehat f_{(\eta ,w_0 ,v_0 )} (\zeta )} d\eta } \right|} ^2 d\zeta
\nn
\en

\be
&&~~\le 4\int\limits_{ - \alpha }^\alpha  {\left| {\int\limits_y^1 {e^{(\eta  - y)\left| \zeta  \right|} \left| {\widehat f_{(\eta ,w_0 ,v_0 )} (\zeta ) - \widehat f_{(\eta ,w^{(\alpha )} ,v_\varepsilon  )} (\zeta )} \right|d\eta } } \right|^2 d\zeta }
\nn\\
&&~~~~+ 4\int\limits_{\left| \zeta  \right| > \alpha }^{} {\left| {\int\limits_y^1 {e^{(\eta  - y)\left| \zeta  \right|} \widehat f_{(\eta ,w_0 ,v_0 )} (\zeta )} d\eta } \right|} ^2 d\zeta
\nn\\
&&~~\le 4\int\limits_{ - \alpha }^\alpha  {\left| {\int\limits_y^1 {e^{(\eta  - y)\alpha } \left| {\widehat f_{(\eta ,w_0 ,v_0 )} (\zeta ) - \widehat f_{(\eta ,w^{(\alpha )} ,v_\varepsilon  )} (\zeta )} \right|d\eta } } \right|^2 d\zeta }
\nn\\
&&~~~~+ 4e^{ - 2y\alpha } \int\limits_{\left| \zeta  \right| > \alpha }^{} {\left| {\int\limits_y^1 {e^{(\eta  - y)\left| \zeta  \right| + y\alpha } \widehat f_{(\eta ,w_0 ,v_0 )} (\zeta )} d\eta } \right|} ^2 d\zeta
\nn\\
&&~~\le 4e^{ - 2y\alpha } \left\{ {\int\limits_{ - \alpha }^\alpha  {\int\limits_y^1 {e^{2\eta \alpha } \left| {\widehat f_{(\eta ,w_0 ,v_0 )} (\zeta ) - \widehat f_{(\eta ,w^{(\alpha )} ,v_\varepsilon  )} (\zeta )} \right|^2 d\eta } d\zeta } } \right.
\nn\\
&&~~~~+ \left. {\int\limits_{\left| \zeta  \right| > \alpha }^{} {\left| {\int\limits_y^1 {e^{(\eta  - y)\left| \zeta  \right| + y\alpha } \widehat f_{(\eta ,w_0 ,v_0 )} (\zeta )} d\eta } \right|} ^2 d\zeta } \right\}
\nn\\
&&~~\le 4e^{ - 2y\alpha } \left\{ {\int\limits_y^1 {e^{2\eta \alpha } \left( {\int\limits_{ - \infty }^{ + \infty } {\left| {f_{(\eta ,w_0 ,v_0 )} (\xi ) - f_{(\eta ,w^{(\alpha )} ,v_\varepsilon  )} (\xi )} \right|^2 d\xi } } \right)d\eta } } \right.
\nn\\
&&~~~~\left. { + \int\limits_{\left| \zeta  \right| > \alpha }^{} {\left| {\int\limits_y^1 {e^{\left| \zeta  \right| + \alpha } \left| {\widehat f_{(\eta ,w_0 ,v_0 )} (\zeta )} \right|} d\eta } \right|} ^2 d\zeta } \right\}
\nn\\
&&~~\le 4e^{ - 2y\alpha } \left\{ {\int\limits_y^1 {e^{2\eta \alpha } \int\limits_{ - \infty }^{ + \infty } {k^2 \left| {v_0 (\xi ,\eta ) - v_\varepsilon  (\xi ,\eta ) + w_0 (\xi ,\eta ) - w^{(\alpha )} (\xi ,\eta )} \right|^2 d\xi } d\eta } } \right.
\nn\\
&&~~~~\left. { + \int\limits_{\left| \zeta  \right| > \alpha }^{} {e^{2(\left| \zeta  \right| + \alpha )} \int\limits_y^1 {\left| {\widehat f_{(\eta ,w_0 ,v_0 )} (\zeta )} \right|^2 d\eta d\zeta } } } \right\}
\nn\\
&&~~\le 4e^{ - 2y\alpha } \left\{ {2k^2 \int\limits_y^1 {e^{2\eta \alpha } \left\| {w_0 (.,\eta ) - w^{(\alpha )} (.,\eta )} \right\|_{L^2 (\RR)}^2 d\eta } } \right.
\nn\\
&&~~~~+ 2\int\limits_0^1 {\int\limits_{ - \infty }^{ + \infty } {e^{2\alpha } k^2 \left| {v_0 (\xi ,\eta ) - v_\varepsilon  (\xi ,\eta )} \right|^2 d\xi d\eta } }
\nn
\en
\be
&&~~~~\left. { + \int\limits_{\left| \zeta  \right| > \alpha }^{} {e^{2(\left| \zeta  \right| + \alpha )} \int\limits_y^1 {\left| {\widehat f_{(\eta ,w_0 ,v_0 )} (\zeta )} \right|^2 d\eta d\zeta } } } \right\}
\nn\\
&&~~\le 4e^{ - 2y\alpha } \left\{ {2k^2 \int\limits_y^1 {e^{2\eta \alpha } \left\| {w_0 (.,\eta ) - w^{(\alpha )} (.,\eta )} \right\|_{L^2 (\RR)}^2 d\eta } } \right.
\nn\\
&&~~~~+ 2k^2 e^{2\alpha } \left\| {v_0  - v_\varepsilon  } \right\|_2^2
\nn\\
&&~~~~\left. { + \int\limits_{\left| \zeta  \right| > \alpha }^{} {e^{2(\left| \zeta  \right| + \alpha )} \int\limits_y^1 {\left| {\widehat f_{(\eta ,w_0 ,v_0 )} (\zeta )} \right|^2 d\eta d\zeta } } } \right\}.
\nn
\en

We put
\be
M_1  = 2k^2 e^{2\alpha } \left\| {v_0  - v_\varepsilon  } \right\|_2^2
\nn
\en
and
\be
M_2  = \int\limits_{\left| \zeta  \right| > \alpha }^{} {e^{2(\left| \zeta  \right| + \alpha )} \int\limits_y^1 {\left| {\widehat f_{(\eta ,w_0 ,v_0 )} (\zeta )} \right|^2 d\eta d\zeta } }.
\nn
\en

Therefore
\be
e^{2y\alpha } \left\| {w_{o(y)}  - w_{}^{(\alpha )} (.,y)} \right\|_{L^2 (\RR)}^2  \le 2k^2 \int\limits_y^1 {e^{2\eta \alpha } \left\| {w_o (.,\eta ) - w^{(\alpha )} (.,\eta )} \right\|_{L^2 (\RR)}^2 d\eta }  + M_1  + M_2.
\nn
\en

Using Gronwall's inequality, we have \be e^{2y\alpha } \left\| {w_o
(.,y) - w^{(\alpha )} (.,y)} \right\|_{L^2 (\RR)}^2  \le \left( {M_1
+ M_2 } \right)e^{2k^2 (1 - y)}, \nn \en hence \be \left\| {w_o
(.,y) - w^{(\alpha )} (.,y)} \right\|_{L^2 (\RR)}^2  \le e^{ -
2\alpha y} (M_1  + M_2 )e^{2k^2 (1 - y)}.
\label{eq054}
 \en

From Theorem 2, we get
\be
e^{ - 2\alpha y} M_1  &=& 2e^{2(1 - y)\alpha } k^2 \left\| {v_0  - v_\varepsilon  } \right\|_2^2
\nn\\
&<& 2D^2 e^{2(1 - y)\alpha } k^2 \left( {\ln \frac{1}{\varepsilon }}
\right)^{ - 2}.
\label{eq055}
 \en

We have
\be
e^{ - 2\alpha y} M_2  &=& e^{ - 2\alpha y} \int\limits_{\left| \zeta  \right| > \alpha }^{} {e^{2(\left| \zeta  \right| + \alpha )} \int\limits_y^1 {\left| {\widehat f_{(\eta ,v_0 ,w_0 )} (\zeta )} \right|^2 d\eta d\zeta } }
\nn
\en
\be
&\le& e^{ - 2\alpha y} \int\limits_{\left| \zeta  \right| > \alpha }^{} {\int\limits_0^1 {\frac{{e^{6\left| \zeta  \right|} }}{{e^{2\alpha } }}\left| {\widehat f_{(\eta ,v_0 ,w_0 )} (\zeta )} \right|^2 } d\eta } d\zeta
\nn\\
&\le& e^{ - 2\alpha (y + 1)} \left\| {e^{3\left| \zeta  \right|}
\left| {\widehat f_{(\eta ,v_0 ,w_0 )} (\zeta )} \right|}
\right\|_2^2.
\label{eq056} \en

From $(\ref{eq054})$-$(\ref{eq056})$ and choosing $\alpha  =
\frac{1}{{2(1 - y)}}\ln \left( {\ln \frac{1}{\varepsilon }}
\right)$, we get \be \left\| {w_o (.,y) - w^{(\alpha )} (.,y)}
\right\|_{L^2 (\RR)}^2  &<& e^{2k^2 } \left[ {2k^2 D^2 \left( {\ln
\frac{1}{\varepsilon }} \right)^{ - 1}  + \left\| {e^{3\left| \zeta
\right|} \left| {\widehat f_{(\eta ,v_0 ,w_0 )} (\zeta )} \right|}
\right\|_2^2 \left( {\ln \frac{1}{\varepsilon }} \right)^{ - 1} }
\right]
\nn\\
&<& C^2 \left( {\ln \frac{1}{\varepsilon }} \right)^{ - 1}.
\nn
\en

Therefore
\be
\left\| {w_o  - w^{(\alpha )} } \right\|_2^{}  < C\left( {\ln \frac{1}{\varepsilon }} \right)^{ - 1/2}.
\nn
\en

Denoting $w_\varepsilon   = w^{(\alpha )}$ and $u_\varepsilon   = v_\varepsilon   + w_\varepsilon$, we get
\be
\left\| {u_0  - u_\varepsilon  } \right\|_2  &=& \left\| {w_0  + v_0  - v_\varepsilon   - w_\varepsilon  } \right\|_2
\nn\\
&\le& \left\| {w_0  - w^{(\alpha )} } \right\|_2  + \left\| {v_0  - v_\varepsilon  } \right\|_2
\nn\\
&\le& C\left( {\ln \frac{1}{\varepsilon }} \right)^{ - 1/2}  + D\left( {\ln \frac{1}{\varepsilon }} \right)^{ - 1}
\nn\\
&\le& (C + D)\left( {\ln \frac{1}{\varepsilon }} \right)^{ - 1/2}.
\nn
\en

This completes the proof of Theorem 3. \ \par \ \par \bf
\section {NUMERICAL RESULTS} \par \ \par \rm
Let the problem
 \be
 &&\Delta u=f(x,y,u),\; x\in  \RR  ,\; y>1 \nn \\
&& u(x,1)=\varphi (x)
\label{eq057}
\en
with
\be
 && f(x,y,u)={\rm arctan(}{\left\vert u\right\vert }
+x^{2}+y^{2}{\rm )}+g(x,y) \nn\\
 && g(x,y)={{4Cx}\over{(x^{2}+4)^{3} }}
+{{-2C}\over{(x^{2}+4)^{2} }} -{\rm arctan(}{{C}\over{x^{2}+4}}
+x^{2}+y^{2}{\rm )},\; C={{4}\over{{\sqrt{ 2\pi  }}  }}\nn \\
&& \varphi (x)={{C}\over{4+x^{2} }}
 \nn \en
The exact solution of $(\ref{eq057})$ is $ u(x,y)={{\displaystyle
C}\over{\displaystyle x^{2}+4}} ,\ C=4/{\sqrt{\displaystyle 2\pi  }}
$ and we solve numerically this problem for $ x\in  \RR  ,\  y>1$ by
using the iterative sequence defined by \ \par \be
u^{(n+1)}(x,y)=-{\int _{-\infty  }  ^{\infty  }  {{\int _{1}
^{\infty  } {N(x,y,}}  \xi  ,\eta  )f(\xi  ,\eta  ,u^{(n)}(\xi ,\eta
))d\xi  d\eta +}}
\nn \\
 {\int _{-\infty  }  ^{\infty  }{N_{y}(x,y,\xi ,1)\varphi  (\xi  )d\xi }}
 \nn \en
 and whose graphic is displayed in Fig.1 on the
interval $ [-4,4]\times  [1,4]$.\ \par

\centerline{
\includegraphics[width=3in]{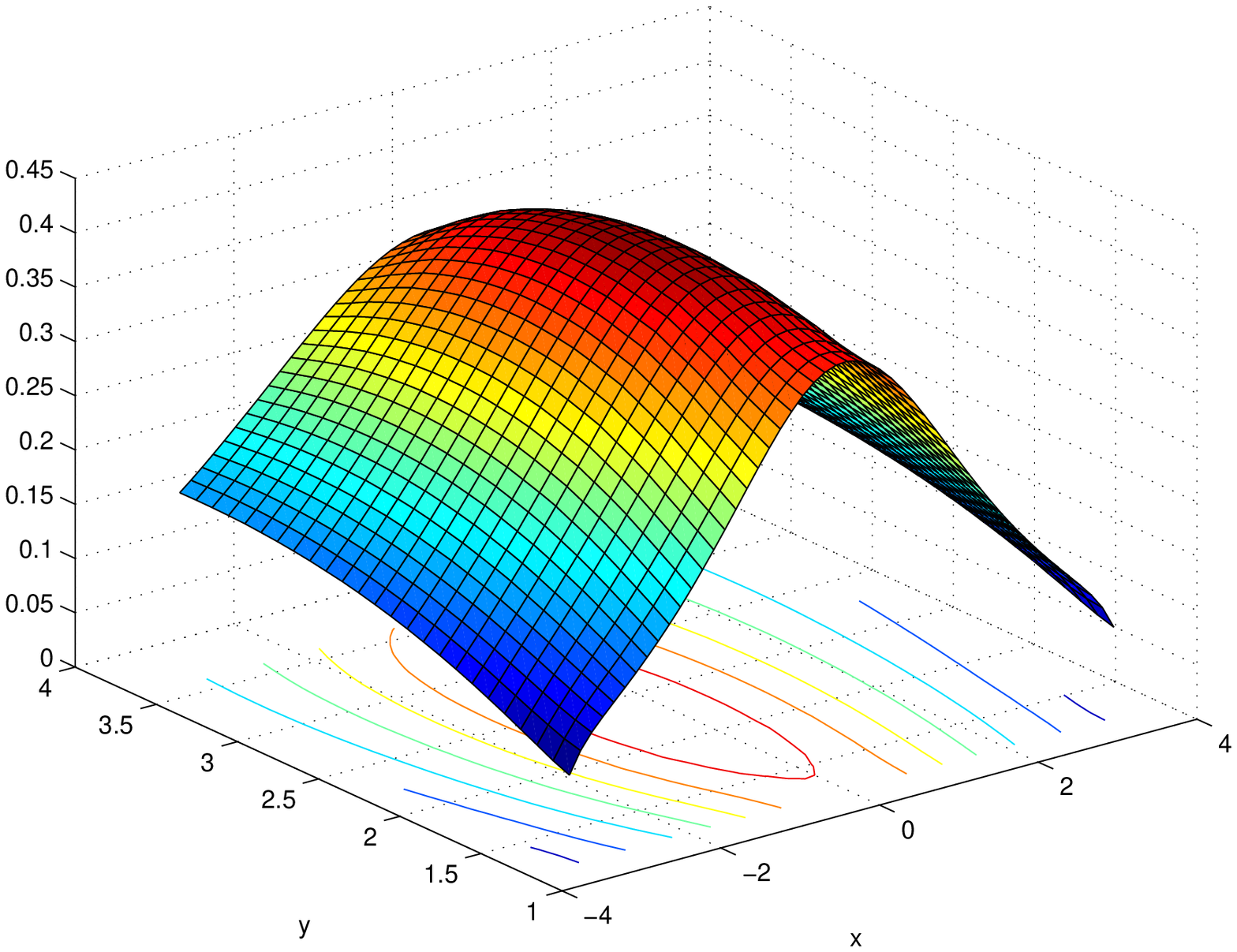}}
\centerline{Fig.1: approximated solution}\par

 For the second problem
 \be
 && \Delta v=0,\; x\in  \RR  ,\; 0<y<1
\nn \\
&& v(x,1)=\varphi (x),\; v_{y}(x,1)=0
 \label{eq058} \en
 whose exact solution is $ v(x,y)={{\displaystyle 1}\over{\displaystyle
2{\sqrt{\displaystyle 2\pi }}  }}  \Big [{{\displaystyle
3-y}\over{\displaystyle x^{2}+(3-y)^{2}}} +{{\displaystyle
1+y}\over{\displaystyle x^{2}+(1+y)^{2}}}  \Big ]$, we calculate the
regularized solution $ v_{\varepsilon  }(x,y)$ of $(\ref{eq058})$
for $ \varepsilon =10^{-2}$ by the formula $(\ref{eq051})$.\ \par
Finally for the third problem \be && \Delta w=f(x,y,v+w),\; x\in
\RR,\;0<y<1
\nn \\
&& w(x,1)=0,\; w_{y}(x,1)=0
 \label{eq059}
\en the regularized solution $ w_{\varepsilon  }(x,y)$ of the
problem $(\ref{eq059})$ is calculated from its definition
$(\ref{eq052})$. So in Fig.2 we have drawed the regularized solution
of the problem 2 i.e. $ u_{\varepsilon }(x,y)=v_{\varepsilon
}(x,y)+w_{\varepsilon }(x,y)$ on the interval $ [-4,4]\times
[0,1]$.\
\par \centerline{
\includegraphics[width=3in]{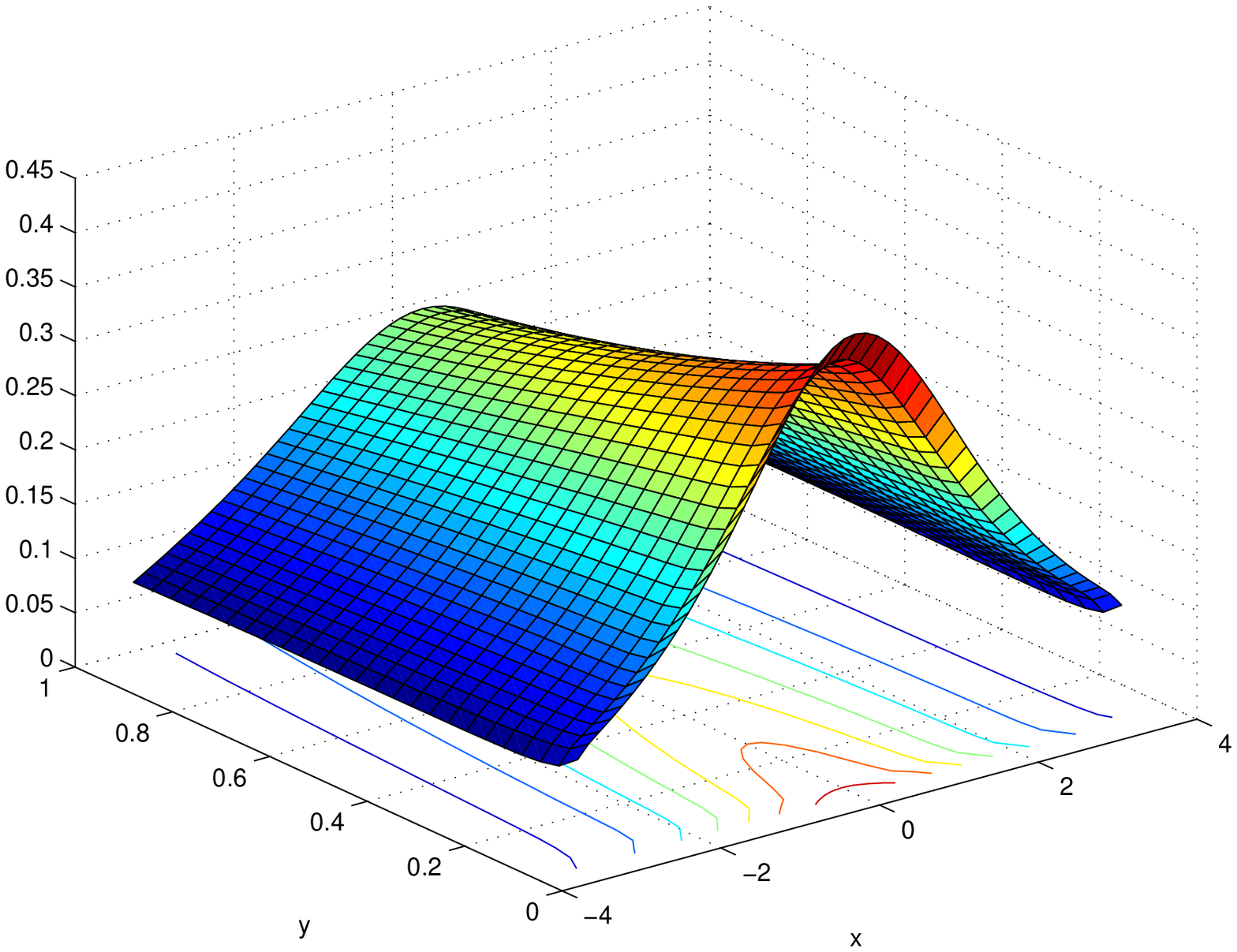}}
\centerline{Fig.2: regularized solution}\par

 For comparison Fig.3 gives the exact solution on the
interval $ [-4,4]\times  [0,4]$. \ \par \centerline{
\includegraphics[width=3in]{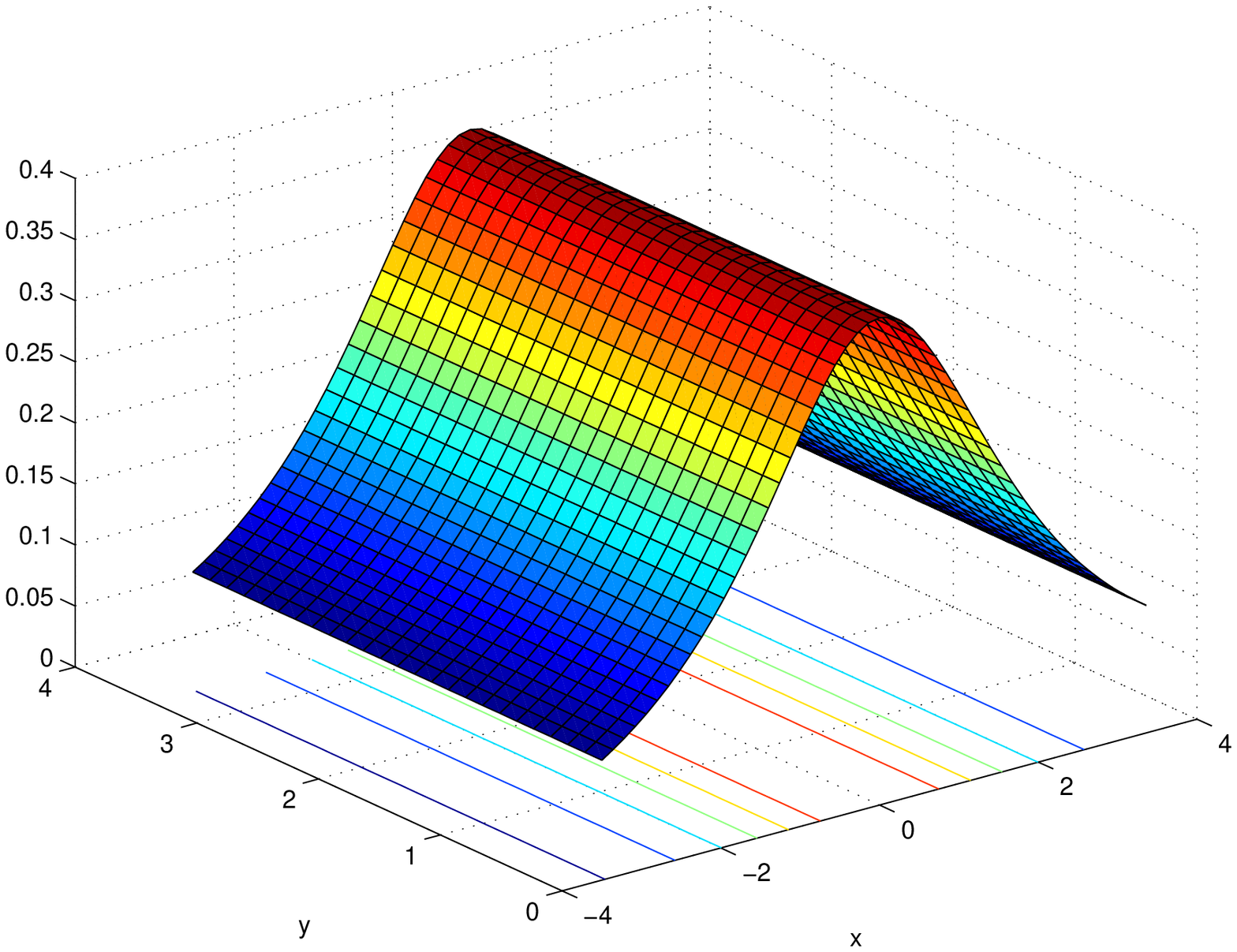}}
\centerline{Fig.3: exact solution}\par

\begin{center}REFERENCES
\end{center}

[AS] Andersen, R. S. and Saull, V. A., {\it Surface temperature history determination from bore hole measurements}, Journal of International Association Mathematical Geology, {\bf 5}, pp. 269-283, 1973.

[B] Bourgeois, L., {\it A mixed formulation of quasi-reversibility to solve the Cauchy problems for Laplace's equations}, Inverse Problems {\bf 21}, No. 3, pp. 1087-1104, 2005.

[BBC] Beck, J. V., Blackwell, B.  and Clair, C. R. St., {\it Inverse Heat Conduction}, Ill-posed Problem, Wiley, New York, 1985.

[C] Colton, D., {\it Partial differential equations}, Random House, New York, 1988.

[C1] Cavazzoni, R., {\it On the Cauchy problem for elliptic equation in a disk}, Rend. Circ. Mat. Palermo (2) {\bf 52}, No. 1, pp. 131-140, 2003.

[CHWY] Cheng, J., Hon, Y. C., Wei, T. and Yamamoto, M., {\it Numerical computation of a Cauchy problem for Laplace's equations}, ZAMM Z. Angew. Math. Mech. {\bf 81}, No. 10, pp. 665-674, 2001.

[HQ] Huang, Y. and Quan, Z., {\it Regularization for a class of ill-posed Cauchy problems}, Proc. Amer. Math. Soc., Vol. 133, No. 10,pp. 3005-3012, 2005.

[HR] Han, H. and Reinhardt, H. J., {\it Some stability estimates for Cauchy problems for elliptic equations}, J. Inverse Ill-Posed Probl. {\bf 5}, No. 5, pp. 437-454, 1997.

[HR1] Hao, D. N. and Reinhardt, H. J., {\it On a sideways parabolic equation}, Inverse Problems, {\bf 13}, pp. 297-309, 1997.

[HRS] Hao, D. N., Reinhardt, H. J. and Schneider, A., {\it Numerical solution to a sideways parabolic equation}, International Journal for Numerical methods in Engineering, {\bf 50}, pp. 1253-1267, 2001.

[KS] Klibanov, M. V. and Santosa, F. V.,{\it A computational quasi-reversibility method for Cauchy problems for Laplace's equation},
SIAM J. Appl. Math. {\bf 51}, No. 6, pp. 1653-1675, 1991.

[K] Krylov, A. L., {\it The Cauchy problem for Laplace's equations in the complex domain}, Dokl. Akad. Nauk  SSSR {\bf 188}, pp. 748-751, 1969.

[L] Levine, H. A., {\it Continuous data dependence, regularization and a three lines theorem for the heat equation with data in a space like direction}, Annali di Matematica Pura ed Applicata, {\bf CXXXIV}, pp. 267-286, 1983.

[LV] Levine H. A. and Vessella S., {\it Estimates and regularization for solutions of some ill-posed problems of elliptic and parabolic type}, Rediconti del Circolo Matematica di Palermo, {\bf 34}, pp. 141-160, 1985.

[M] Murio, D., {\it The mollification method and the numerical solution of ill-posed problems}, Wiley, New York, 1993.

[QD] P. H. Quan and N. Dung, {\it A backward nonlinear heat equation: regularization with error estimates}, accepted for publication in Applicable Analysis, 2005.

[QT] P. H. Quan and D. D. Trong, {\it Temperature Determination from Interior Measurements: the Case of Temperature Nonlinearly Dependent Heat Source}, Vietnam Journal of Math., {\bf 32}, pp. 131-142, 2004.

[T] Tautenhahn, U., {\it Optimal stable solution of Cauchy problems for elliptic equations}, Z. Anal. Anwendungen {\bf 15}, No. 4, pp. 961-984, 1996.

[TT] D. D. Trong and N. H. Tuan, {\it Regularization and error estimates for nonhomogeneous backward heat problems}, Electronic J. of Diff. Eq. Vol. {\bf 2006}, No. 04, pp. 1-10, 2006.

\end{document}